\documentclass{llncs}

\usepackage{epsf,amssymb,amsmath}

\def\bydef{:=}
\def\F{\mathcal{F}}
\def\G{\mathcal{G}}
\def\H{\mathcal{H}}

\def\R{\mathbb{R}}

\def\eps{\varepsilon}

\newcommand{\e}[1]{\mathbb{E}\left[#1 \right]}

\newcommand{\ee}[2]{\mathbb{E}_{#1}\left[#2 \right]}

\newcommand{\abs}[1]{\left|#1\right|}
\newcommand{\norm}[1]{\left\|#1\right\|}

\newcommand{\inner}[1]{\left\langle#1\right\rangle}

\DeclareMathOperator{\conv}{conv}

\title{Some Local Measures of Complexity of Convex Hulls and Generalization
Bounds}
\author{Olivier Bousquet\inst{1} \and Vladimir Koltchinskii\inst{2}\thanks{Partially supported by NSA Grant MDA904-99-1-0031} \and
Dmitriy Panchenko\inst{2}}
\institute{
Centre de Math\'ematiques Appliqu\'ees\\
Ecole Polytechnique\\
91128 Palaiseau, FRANCE\\
\email{bousquet@cmapx.polytechnique.fr}
\and Department of Mathematics and Statistics\\
The University of New Mexico\\
Albuquerque, NM 87131-1141, U.S.A.\\
\email{\{vlad,panchenk\}@math.unm.edu}}

\begin{document}
\maketitle
\begin{abstract}
We investigate measures of complexity of function classes based on
continuity moduli of Gaussian and Rademacher processes.
For Gaussian processes, we obtain bounds on the continuity modulus on the
convex hull of a function class in terms of the same quantity for the class
itself.
We also obtain new bounds on generalization error in terms of localized
Rademacher complexities. This allows us to prove new results about
generalization performance for convex hulls in terms of characteristics of
the base class.
As a byproduct, we obtain a simple proof of some of the known bounds on the
entropy of convex hulls.
\end{abstract}

\section{Introduction}
Convex hulls of function classes have become of great interest in Machine
Learning since the introduction of AdaBoost and other methods of combining
classifiers.
The most commonly used measure of complexity of convex hulls is based on
covering numbers (or metric entropies). The first bound on the entropy of the
convex hull of a set in a Hilbert space was obtained by Dudley
\cite{dudley87} and later refined by Ball and Pajor \cite{bp90} and
a different proof was given independently by van der Vaart and Wellner
\cite{vdVW96}.
These authors considered the case of polynomial growth of the covering
numbers of the base class. Sharp bounds in the case of exponential growth of
the covering numbers of the base class as well as extension of previously
konwn results to the case of Banach spaces were obtained later
\cite{carl01,mendelson01,lilinde01,gao01,cre01}.

In Machine Learning, however, the quantities of primary importance for
determining the generalization performance are not the entropies themselves
but rather localized Gaussian or Rademacher complexities of the function
classes \cite{kolpan00,barboumen02}.
These quantities are closely related to continuity moduli of the
corresponding stochastic processes.

Our main purpose in this paper is to provide an easy bound on the continuity
modulus of stochastic processes like Rademacher or
Gaussian processes on the convex hull of a class in terms of the continuity
modulus on the class itself.
We combine this result with some new bounds on the generalization error in
function learning problems based on localized Rademacher complexities.
This allows us to bound the generalization error for convex hulls in terms
of characteristics of the base class.

In addition to this, we use the bounds on continuity moduli on convex hulls
to give very simple proofs of some previously known results on the entropy of
such classes.

\section{Continuity Modulus on Convex Hulls}\label{sec:main}
Let $\F$ be a subset of a Hilbert space $\H$ and $W$ denote an isonormal
Gaussian process defined on $\H$, that is a collection $(W(h))_{h\in\H}$ of
Gaussian random variables indexed by $\H$ such that 
\[
\forall h\in\H,\,\e{W(h)} = 0\,\, \mbox{ and }\,\, \forall h,h'\in\H,\,
\e{W(h)W(h')} = \inner{h,h'}_\H\,.
\]
We define the modulus of continuity of the process $W$ as
\[
\omega(\F,\delta)\bydef \omega_\H(\F,\delta) = \e{\sup_{f,g\in\F \atop
\|f-g\|\leq \delta} \abs{W(f)-W(g)}}\,.
\]
Let $\F_\eps$ denote a minimal $\eps$-net of $\F$, i.e. a subset of $\F$ of
minimal cardinality such that $\F$ is contained in the union of the balls of
radius $\eps$ with centers in $\F_\eps$.
Let $\F^\eps$ denote a maximal $\eps$-separated subset of $\F$, i.e. a
subset of $\F$ of maximal cardinality such that the distance between any two
points in this subset is larger than or equal to $\eps$.
The $\eps$-covering number of $\F$ is then defined as
\[
N(\F,\eps)\bydef N_\H(\F,\eps) = |\F_\eps|\,,
\]
and the $\eps$-entropy is $H(\F,\eps)=\log N(\F,\eps)$.

\subsection{Main Result}
Our main result relates the continuity modulus on the convex hull of a set
$\F$ to the continuity modulus on this set.
\begin{theorem}\label{th:mod}
We have for all $\delta\geq 0$
\[
\omega(\conv(\F),\delta) \leq \inf_{\eps}
\left(2\omega(\F,\eps) + \delta \sqrt{N(\F,\eps)}\right)\,.
\]
\end{theorem}
\begin{proof}
Let $\eps>0$, $L$ be the linear span of $\F_\eps$ and
$\Pi_L$ be the orthogonal projection on $L$. We have for all $f\in\F$,
\[
f = \Pi_L(f) + \Pi_{L^\perp}(f)\,.
\]
\begin{eqnarray*}
\omega(\conv(\F),\delta) &\leq& \e{\sup_{f,g\in\conv(\F)\atop \|f-g\|\leq \delta}
\abs{W(\Pi_L f)-W(\Pi_L g)}}\\ 
&&+ \e{\sup_{f,g\in\conv(\F) \atop \|f-g\|\leq \delta}
\abs{W(\Pi_{L^\perp}f)-W(\Pi_{L^\perp}g)}}\,.
\end{eqnarray*}
Now since for any orthogonal projection $\Pi$, $\norm{\Pi(f)-\Pi(g)}\leq
\norm{f-g}$ we have
\[
\omega(\conv(\F),\delta) \leq \omega(\Pi_L \conv(\F),\delta) +
\omega(\Pi_{L^\perp} \conv(\F),\delta)\,.
\]
Moreover, we have $\Pi \conv(\F) = \conv(\Pi \F)$ by linearity of the
orthogonal projection so that
\[
\omega(\conv(\F),\delta) \leq \omega(\conv(\Pi_L \F),\delta) +
\omega(\conv(\Pi_{L^\perp}\F),\delta)\,.
\]
This gives the first inequality.
Next we have
\[
\omega(\Pi_L \conv(\F),\delta) \leq \omega(L,\delta)\,,
\]
and by linearity of $W$,
\[
\omega(L,\delta)=\e{\sup_{f\in L\atop \norm{f}\leq \delta}\abs{W(f)}} \leq
\delta \e{\sup_{\norm{y}_{\R^d}\leq 1\atop y\in \R^d}\inner{Z,y}}\,,
\]
where $Z$ is a standard normal vector in $\R^d$ (with $d=\dim L$ and
$\norm{\cdot}_{\R^d}$ the euclidean norm in $\R^d$).
This gives
\[
\omega(L,\delta)\leq \delta\e{\norm{Z}_{\R^d}} \leq \delta \sqrt{\dim L} \leq
\delta\sqrt{N(\F,\eps)}\,.
\]
We also get
\[
\omega(\Pi_{L^\perp} \conv(\F),\delta) \leq
2\e{\sup_{f\in\conv(\F)}\abs{W(\Pi_{L^\perp}f)}}\,.
\]
Since $\Pi_{L^\perp}$ is linear, the supremum is attained at elements of
$\F$, that is
\[
\omega(\Pi_{L^\perp} \conv(\F),\delta) \leq
2\e{\sup_{f\in \F}\abs{W(\Pi_{L^\perp}f)}}\,.
\]
Now for each $f\in\F$, let $g$ be the closest point to $f$ in
$\F_\eps$.
Then we have $\norm{f-g}\leq \eps$ and $g\in L\cap \F$ so that
$\Pi_{L^\perp}g=0$ and thus
\[
\omega(\Pi_{L^\perp} \conv(\F),\delta) \leq
2\e{\sup_{f,g\in \F \atop \norm{f-g}\leq
\eps}\abs{W(\Pi_{L^\perp}f)-W(\Pi_{L^\perp}g)}}\,.
\]
Now since $\Pi_{L^\perp}$ is a contraction, using Slepian's lemma (see
\cite{lt91}, Theorem 3.15 page 78) we get
\[
\omega(\Pi_{L^\perp} \conv(\F),\delta) \leq
2\e{\sup_{f,g\in \F \atop \norm{f-g}\leq \eps}\abs{W(f)-W(g)}} = 2\omega(\F,\eps)\,.
\]
This concludes the proof.\qed
\end{proof}

Note that Theorem \ref{th:mod} allows us to give a positive answer to a
question raised by Dudley \cite{dud}. Indeed, we can prove that the
convex hull of a uniformly Donsker class is uniformly Donsker. Due to
lack of space, we do not give the details here.

\subsection{Examples}
As an application of Theorem \ref{th:mod}, we will derive bounds on the
continuity modulus of convex hulls of classes for which we know the rate of
growth of the entropy.

By Dudley's entropy bound (see \cite{lt91}, Theorem 11.17, page 321) we have
\[
\omega(\F,\eps) \leq K\int_0^\eps H^{1/2}(\F,u)\,du\,.
\]
We will also use below the following version of this result (that easily
follows from Dudley's chaining argument and is well known)
\[
\omega(\F^{\delta},\eps) \leq K\int_\delta^\eps H^{1/2}(\F^{\delta},u)\,du\,,
\]
for all $\eps>\delta$.

We first consider the case when the entropy of the base class grows
logarithmically.
\begin{example}\label{ex:e1}
If for all $\eps>0$,
\[
N(\F,\eps) \leq K\eps^{-V}\,,
\]
then for all $\delta>0$,
\[
\omega(\conv(\F), \delta) \leq K\delta^{2/(2+V)}\log^{V/(2+V)}\delta^{-1}\,.
\]
\end{example}
\begin{proof}
We have from Theorem \ref{th:mod},
\begin{eqnarray*}
\omega(\conv(\F),\delta) &\leq& \inf_\eps\left(K\int_0^\eps \log^{1/2} u^{-1}
du + \delta\eps^{-V/2}\right)\\
&\leq& \inf_\eps\left(K\eps\log^{1/2}\eps^{-1} + \delta\eps^{-V/2}\right)\,.
\end{eqnarray*}
Choosing
\[
\eps = \delta^{2V/(2+V)}\log^{2V/(2+V)} \delta^{-1}\,,
\]
we obtain for $\delta\leq 1$,
\[
\omega(\conv(\F),\delta) \leq K\delta^{2/(2+V)}\log^{V/(2+V)} \delta^{-1}\,.
\]\qed
\end{proof}
Although the main term in the above bound is correct, we obtain a
superfluous logarithm. This logarithm can be removed if one uses
directly the entropy integral in combination with results on the
entropy of the convex hull of such classes
\cite{bp90,vdVW96,mendelson01}. At the moment of this writing, we do
not know a simple proof of this fact that does not rely upon the bounds on
the entropy of convex hulls.

Now we consider the case when the entropy of the base class has polynomial
growth.
In this case, we shall distinguish several situations: when the exponent is
larger than $2$, the class is no longer pre-Gaussian which means that the
continuity modulus is unbounded. However, it is possible to study the
continuity modulus of a restricted class. Here we consider the convex hull
of a $\delta$-separated subset of the base class, for which the continuity
modulus is bounded when computed at a scale proportional to $\delta$.
\begin{example}\label{ex:e2}
If for all $\eps>0$,
\[
H(\F,\eps) \leq K\eps^{-V}\,,
\]
then for all $\delta>0$, for $0<V<2$,
\[
\omega(\conv(\F), \delta) \leq K\log^{1/2-1/V}\delta^{-1}\,,
\]
for $V=2$,
\[
\omega(\conv(\F^{\delta/4}), \delta) \leq K\log \delta^{-1}\,,
\]
and for $V>2$,
\[
\omega(\conv(\F^{\delta/4}), \delta) \leq K\delta^{1-V/2}\,.
\]
\end{example}
\begin{proof}
We have from Theorem \ref{th:mod}, for $\epsilon>\delta/4$,
\[
\omega(\conv(\F^{\delta/4}),\delta) \leq
\inf_\eps\left(K\int_{\delta/4}^\eps u^{-V/2} du +
\delta\exp(K\eps^{-V}/2)\right)\,.
\]
For $0<V<2$, this gives
\[
\omega(\conv(\F),\delta) \leq \inf_\eps\left(K\eps^{(2-V)/2} +
\delta\exp(K\eps^{-V}/2)\right)\,.
\]
Choosing
\[
\eps = K^{1/V}\log^{-1/V} \delta^{-1}\,,
\]
we obtain for $\delta$ small enough
\[
\omega(\conv(\F),\delta) \leq K\log^{(V-2)/2V} \delta^{-1}\,.
\]
For $V=2$, we get
\[
\omega(\conv(\F^{\delta/4}),\delta) \leq \inf_\eps\left(K\log
\frac{4\epsilon}{\delta} + \delta\exp(K\eps^{-2}/2)\right)\,.
\]
Taking $\epsilon=1/4$ we get for $\delta$ small enough
\[
\omega(\conv(\F^{\delta/4}),\delta) \leq K\log \delta^{-1}\,.
\]
For $V>2$, we get
\[
\omega(\conv(\F^{\delta/4}),\delta) \leq
\inf_\eps\left(K\delta^{(2-V)/2} - \eps^{(2-V)/2} +
\delta\exp(K\eps^{-2}/2)\right)\,.
\]
Taking $\eps\rightarrow \infty$, we obtain
\[
\omega(\conv(\F^{\delta/4}),\delta) \leq K\delta^{(2-V)/2}\,.
\]\qed
\end{proof}

\section{Generalization Error Bounds}

\subsection{Results}
We begin this section with a general bound that relates the error of the
function minimizing the empirical risk to a local measure of complexity of
the class which is the same in spirit as the bound in \cite{kolpan00}.

Let $(S,\mathcal{A})$ be a measurable space and let $X_1,\ldots,X_n$ be $n$
i.i.d. random variables in this space with common distribution $P$.
$P_n$ will denote the empirical measure based on the sample
\[
P_n = \frac{1}{n}\sum_{i=1}^n \delta_{X_i}\,.
\]
In what follows, $\H=L_2(P_n)$ and we are using the notations of Section
\ref{sec:main}.

We consider a class $\F$ of measurable functions defined on $S$ with values
in $[0,1]$.
We assume in what follows that $\F$ also satisfies standard measurability
conditions used in the theory of empirical processes as in
\cite{dudley00,vdVW96}.

We define
\[
R_n(f)\bydef \frac{1}{n} \sum_{i=1}^n \eps_i f(X_i)\,,
\]
and let $\psi_n$ be an increasing concave (possibly data-dependent random)
function with $\psi_n(0)=0$ such that 
\[
\ee{\eps}{\sup_{P_n f\leq r}\abs{R_n(f)}}\leq
\psi_n(\sqrt{r}),\,\, \forall r\geq 0\,.
\]
Let $\hat{r}_n$ be the largest solution of the equation
\begin{equation}\label{Uo}
r = \psi_n(\sqrt{r})\,.
\end{equation}
The solution $\hat{r}_n$ of (\ref{Uo}) gives what is usually called 
zero error rate for the class $\F$\cite{kolpan00}, i.e. the bound for $Pf$
given that $P_n f=0$.

The bounds we obtain below are data-dependent and they do not require any
structural assumptions on the class (such as VC conditions or entropy
conditions). Note that $\hat{r}_n$ is determined only by the restriction
of the class $\F$ to the sample $(X_1,\ldots,X_n)$.
\begin{theorem}\label{th:dmi}
If $\psi_n$ is a non-decreasing concave function and $\psi_n(0)=0$
then there exists $K>0$ such that with probability at least
$1-2e^{-t}$ for all $f\in\F$
\begin{equation}\label{Pnf}
Pf \leq K\left(P_n f + \hat{r}_n + \frac{t+\log\log n}{n}\right)\,.
\end{equation}
\end{theorem}
It is most common to estimate the expectation of Rademacher processes 
via entropy integral (Theorem 2.2.4 in \cite{vdVW96}):
\[
\ee{\eps}{\sup_{P_n f\leq \delta}\abs{R_n(f)}}
\leq \frac{4\sqrt{3}}{\sqrt{n}} 
\int_{0}^{\sqrt{\delta}/2}H^{1/2}(\F,u)du\,,
\]
which means one can choose $\psi_n(\delta)$ as the right hand side of
the above bound.
This approach was used for instance in \cite{kolpan00}.

Our goal here will be to apply the bound of Theorem \ref{th:dmi} to the
function learning problem in the convex hull of a given class.

Let $\G$ be a class of measurable functions from $S$ into $[0,1]$.
Let $g_0\in\conv(\G)$ be an unknown target function.
The goal is to learn $g_0$ based on the data
$(X_1,g_0(X_1)),\ldots,(X_n,g_0(X_n))$.
We introduce $\hat{g}_n$ defined as
\[
\hat{g}_n\bydef \arg\min_{g\in\conv(\G)} P_n|g-g_0|\,,
\]
which in principle can be computed from the data.

We introduce the function $\psi_n(\G,\delta)$ defined as
\[
\psi_n(\G,\delta) \bydef \sqrt{\frac{\pi}{2n}}\inf_{\eps>0}
\left(\omega(\G,\eps) + \delta\sqrt{N(\G,\eps)}\right)\,.
\]

\begin{corollary}
Let $\hat{r}_n(\G)$ be the largest solution of the equation
\[
r = \psi_n(\G,\sqrt{r})\,.
\]
Then there exists $K>0$ such that for all $g_0\in\conv(\G)$ the following
inequality holds with probability at least $1-2e^{-t}$
\[
P|\hat{g}_n-g_0| \leq K\left(\hat{r}_n(\G) + \frac{t+\log\log n}{n}\right)\,.
\]
\end{corollary}
\begin{proof}
Let $\F=\{|g-g_0|:g\in\conv(\G)\}$.
Note that $\psi_n(\G,\delta)$ is concave non-decreasing (as the infimum of
linear functions) and $\psi_n(\G,0)=0$, it can thus be used in Theorem
\ref{th:dmi}. 
We obtain (using bound (4.8) on page 97 of \cite{lt91})
\begin{eqnarray*}
\e{\sup_{f\in\F\atop P_nf\leq r} \abs{R_n(f)}}
&\leq& \sqrt{\frac{\pi}{2n}} \e{\sup_{f\in\F\atop P_nf\leq r} \abs{W_{P_n}(f)}}\\ 
&\leq& \sqrt{\frac{\pi}{2n}} \e{\sup_{f\in\F\atop (P_nf^2)^{1/2}\leq \sqrt{r}}
\abs{W_{P_n}(f)}}\\
&\leq& \sqrt{\frac{\pi}{2n}} \omega(\conv{\G},\sqrt{r})
\leq \psi_n(\G,\sqrt{r})\,,
\end{eqnarray*}
where in the last step we used Theorem \ref{th:mod}.
To complete the proof, it is enough to notice that $P_n|\hat{g}_n-g_0|=0$
(since $g_0\in\conv(\G)$) and to use the bound of Theorem \ref{th:dmi}.
\qed
\end{proof}

A simple application of the above corollary in combination with the bounds
of examples \ref{ex:e1} and \ref{ex:e2} give, for instance, the following
rates.
If the covering numbers of the base class grow polynomially, i.e. for
some $V>0$,
\[
N(\G,\eps)\leq K\eps^{-V}\,,
\]
then we obtain $\hat{r}_n$ of the order of
\[
n^{-\frac{1}{2}\frac{2+V}{1+V}}\,.
\]
This can be compared with the main result in \cite{mendelson02}.
If the entropy is polynomial with exponent $0<V<2$, $\hat{r}_n$ is of
the order of
\[
n^{-\frac{1}{2}}\log^{1/2-1/V}n\,.
\]

\subsection{Additional Proofs}
Our main goal in this section is to prove Theorem \ref{th:dmi}.

Denote 
\[
l(\delta)=2\log\left(\frac{\pi}{\sqrt{3}}\log_2 \frac{2}{\delta}\right)
\]
and define $U(\delta)$ as the largest solution of the equation
\begin{equation}\label{U}
U = \delta+8\ee{\eps}{\sup_{P_n f\leq U}\abs{
R_n(f)}}+
\left(\frac{2\delta(t+l(\delta))}{n}\right)^{1/2}+
\frac{10(t+l(\delta))}{3n}
\end{equation}
while $r(\delta)$ is the largest solution of the equation
\begin{equation}\label{r}
r = \delta + 8\ee{\eps}{\sup_{P_n f\leq U(2r)}\abs{R_n(f)}}+
\left(\frac{4r(t+l(2r))}{n}\right)^{1/2}+
\frac{10(t+l(2r))}{3n}.
\end{equation}
Notice that the construction of $r(\delta)$ depends only
on the sample $(X_1,\ldots,X_n)$ and the restriction of the class 
$\F$ to the sample.
\begin{theorem}
With probability at least $1-2e^{-t}$ for all $f\in\F$
\begin{equation}
Pf\leq r(P_n f).
\label{Pf}
\end{equation}
\end{theorem}
\begin{proof}
We define $\delta_k=2^{-k}$ for $k\geq 0,$ and consider a 
sequence of classes 
\[
\F_k=\{f\in \F : \delta_{k+1}<Pf\leq \delta_k\}\,.
\]
If we denote
\[
R_k=\ee{\eps}{\sup_{\F_k}\abs{R_n(f)}}\,,
\]
then the symmetrization inequality implies that
\[
\e{\sup_{\F_k}\abs{P_n f - Pf}} \leq 2\e{R_k}\,,
\]
which in combination with Theorem 3 in \cite{bousquet02} (with
$P(f-Pf)^2\leq Pf^2\leq Pf\leq \delta_k$)
implies that with probability at least $1-e^{-t}$
for all $f\in \F_k$
\[
\abs{P_nf-Pf}\leq 4\e{R_k} + \left(\frac{2\delta_k t}{n}\right)^{1/2}+
\frac{4t}{3n}\,.
\]
Theorem 16 in \cite{boulugmas02} gives that with probability at least
$1-e^{-t}$
\[
\e{R_k} \leq \left(\left(\frac{t}{2n}\right)^{1/2}+
\left(\frac{t}{2n}+R_k\right)^{1/2}\right)^2 \leq
\frac{2t}{n}+2R_k.
\]
Therefore, with probability at least $1-2e^{-t}$ 
for all $f\in \F_k$
\[
\abs{P_nf-P f} \leq 8 R_k +\left(\frac{2\delta_k t}{n}\right)^{1/2}+
\frac{10 t}{3n}.
\]
Finally, replacing $t$ by $t+l(\delta_k)$ and applying the union bound
we get that with probability at least $1-2e^{-t}$
for all $k\geq 0$ and for all $f\in \F_k$
\begin{equation}
\abs{P_nf-P f} \leq 8 R_k +\left(\frac{2\delta_k
(t+l(\delta_k))}{n}\right)^{1/2}+ \frac{10 (t+l(\delta_k))}{3n}.
\label{main}
\end{equation}
If we denote
\[
U_k = \delta_k + 8 R_k +\left(\frac{2\delta_k 
(t+l(\delta_k))}{n}\right)^{1/2}+
\frac{10 (t+l(\delta_k))}{3n}
\]
then on this event for any fixed $k$ and for all $f\in \F_k,$
$P_n f\leq U_k$ and, hence,
\[
R_k \leq \ee{\eps}{\sup_{P_n f\leq U_k} \abs{R_n(f)}}
\]
which can be rewritten in terms of $U_k$ as
\[
U_k \leq \delta_k+ 8\ee{\eps}{\sup_{P_n f\leq U_k} 
\abs{R_n(f)}} + \left(\frac{2\delta_k (t+l(\delta_k))}{n}\right)^{1/2}+
\frac{10 (t+l(\delta_k))}{3n}.
\]
This means that $U_k\leq U(\delta_k),$ where $U(\delta)$
is defined in (\ref{U}).
Finally, (\ref{main}) implies that for all $k$ and $f\in \F_k$ 
\[
Pf \leq P_n f+ 8\ee{\eps}{\sup_{P_n f\leq U(\delta_k)} 
\abs{R_n(f)}} + \left(\frac{2\delta_k (t+l(\delta_k))}{n}\right)^{1/2}+
\frac{10 (t+ l(\delta_k))}{3n}.
\]
If $f\in \F_k$ then $\delta_k\leq 2Pf,$ which proves the theorem.
\qed
\end{proof}

Notice that if we replace the right-hand sides of (\ref{U}) and
(\ref{r}) by upper bounds, we only increase the value of the solutions
and the theorem remains true for these new solutions.
Moreover, since the solution of (\ref{r}) is necessarily larger than
$1/n$, it is enough to consider (\ref{U}) only for $\delta>1/n$.
So assuming that we have the bound
\[
\ee{\eps}{\sup_{P_n f\leq r}\abs{R_n(f)}}
\leq \psi_n(\sqrt{r})\,,
\]
we can replace (using that $2\sqrt{ab}\leq a + b$) (\ref{U}) and (\ref{r}) by
\begin{equation}
U = K_1\left(\delta + \psi_n(\sqrt{U}) + r_0\right)\,,
\label{Uent}
\end{equation}
\begin{equation}
r = \delta + K_2\left(\psi_n(\sqrt{U^{e}(2r)}) + \sqrt{r r_0} + r_0\right)\,.
\label{rent}
\end{equation}
where $r_0=(t+\log\log n)/n$.
The solutions of those equations are denoted respectively $U_1(\delta)$ and
$r_1(\delta)$.

\flushleft{\em Proof of Theorem \ref{th:dmi}}.
Let $\alpha<1$ and consider $k$ non-negative functions $\phi_i$
satisfying one of the following conditions
\begin{equation}\label{eq:c1}
\forall x>0, \, \forall C>1,\, \phi_i(Cx)\leq C^\alpha\phi_i(x)\,,
\end{equation}
or
\begin{equation}\label{eq:c2}
\phi_i(x) \,\mbox{ is non-increasing for }\, x>0\,.
\end{equation}
Define now for each $i=1,\ldots,k$ $u_i$ as the largest solution
of the equation
\[
u = \phi_i(u)\,,
\]
(assuming the existence of the solutions).

Note that from the conditions (\ref{eq:c1}) or (\ref{eq:c2}), we
obtain for all $c>0$ and all $C>1$
\begin{equation}\label{eq:c3}
\phi_i(C(u_i+c)) \leq C^\alpha(u_i+c)\,.
\end{equation}
We thus deduce that the largest solution
$u^*$ of the equation
\[
u = \sum_{i=1}^k \phi_i(u)\,,
\]
satisfies $u^*\leq C\sum_{i=1}^k u_i$ for some large enough
$C$.

It is easy to see that the right-hand side of (\ref{Uent}) is a sum of
functions satisfying (\ref{eq:c3}). Indeed, we have by the concavity
of $\psi_n$ (and $\psi_n(0)=0$) and the definition of $\hat{r}_n$,
\[
\psi_n(\sqrt{C(\hat{r}_n+c)}) \leq  \sqrt{C}
\psi_n(\sqrt{\hat{r}_n+c}) \leq \sqrt{C}(\hat{r}_n+c)\,.
\]
The above reasoning thus proves that
$U_1(\delta)\leq K(\delta+ \hat{r}_n + r_0))$.

We can thus replace equation (\ref{rent}) by the following whose
solution $r_2(\delta)$ will upper bound $r_1(\delta)$:
\[
r = \delta + K_1\left(\psi_n(\sqrt{K_2(r+\hat{r}_n+r_0)})
+\sqrt{r r_0} + r_0\right)\,.
\]
Once again we can check that the righ-hand side is a sum of
functions satisfying (\ref{eq:c3}).
The same reasoning as before proves that
\[
r(\delta)\leq r_2(\delta)\leq K\left(\delta + \hat{r}_n + r_0\right)\,,
\]
which finishes the proof.
\qed

\section{Entropy of Convex Hulls}
\subsection{Relating Entropy With Continuity Modulus}
By Sudakov's minoration (see \cite{lt91}, Theorem 3.18, page 80) we have
\[
\sup_{\eps>0} \eps H^{1/2}(\F,\eps) \leq K\e{\sup_{f \in\F}\abs{W(f)}}\,.
\]
Let $B(f,\delta)$ be the ball centered in $f$ of radius $\delta$.
We define
\[
H(\F,\delta,\eps)\bydef \sup_{f\in\F} H(B(f,\delta)\cap\F, \eps)\,.
\]
The following lemma relates the entropy of $\F$ with the modulus of
continuity of the process $W$. This type of bound is well known (see
e.g. \cite{lif}) but we give the proof for completeness.
\begin{lemma}\label{le:lif}
Assume $\F$ is of diameter $1$.
For all integer $k$ we have
\[
H^{1/2}(\F,2^{-k}) \leq K\sum_{i=0}^k 2^{i} \omega(\F,2^{1-i})\,.
\]
This can also be written
\[
H^{1/2}(\F,\delta) \leq K\int_\delta^1 u^{-2}\omega(\F,u)\,du\,.
\]
\end{lemma}
\begin{proof}
We have
\begin{eqnarray*}
\omega(\F,\delta)
&=& \e{\sup_{f,g\in\F \atop \|f-g\|\leq \delta}
\abs{W(f)-W(g)}}\\
&\geq& \sup_{f\in\F} \e{\sup_{g\in B(f,\delta)\cap\F}
\abs{W(f)-W(g)}}\\
&\geq& \sup_{f\in\F} \sup_{\eps>0} \eps H^{1/2}(B(f,\delta)\cap\F,\eps)\,,
\end{eqnarray*}
so that we obtain
\[
\frac{\delta}{2}H^{1/2}(\F,\delta,\frac{\delta}{2}) \leq K\omega(\F,\delta)\,.
\]
Notice that we can construct a $2^{-k}$ covering of $\F$ by
covering $\F$ by $N(\F,1)$ balls of radius $1$ and then covering the
intersection of each of these balls with $\F$ with $N(B(f,1)\cap
\F,1/2)$ balls of radius $1/2$ and so on.
We thus have
\[
N(\F,2^{-k})\leq \prod_{i=0}^k \sup_{f\in\F} N(B(f,2^{1-i})\cap\F,2^{-i})\,.
\]
Hence
\[
H(\F,2^{-k}) \leq \sum_{i=0}^k H(\F,2^{1-i},2^{-i})\,.
\]
We thus have
\[
H^{1/2}(\F,2^{-k}) \leq \sum_{i=0}^k H^{1/2}(\F,2^{1-i},2^{-i}) \leq
K\sum_{i=0}^k 2^{i} \omega(\F,2^{1-i})\,,
\]
which concludes the proof.\qed
\end{proof}

Next we present a modification of the previous lemma that can be applied to
$\delta$-separated subsets.
\begin{lemma}\label{le:lif2}
Assume $\F$ is of diameter $1$.
For all integer $k$ we have
\[
H^{1/2}(\F,2^{-k}) \leq K\sum_{i=0}^k 2^{i} \omega(\F^{2^{-i-1}},2^{2-i})\,.
\]
\end{lemma}
\begin{proof}
Notice that for $f\in \F$, there exists $f'\in \F^{\delta/4}$
such that
\[
B(f,\delta)\cap\F \subset B(f',\delta+\delta/4)\cap\F \,.
\]
Moreover, since a maximal $\delta$-separated set is a $\delta$-net, 
\[
N(\F,\delta) \leq \abs{N^{\delta}} = N(\F^{\delta},\delta/2)\,,
\]
since for a $\delta$-separated set $A$ we have $N(A,\delta/2)=|A|$.

Let's prove that we have for any $\gamma$,
\[
\abs{(B(f,\gamma)\cup\F)^{\delta/2}} \leq
\abs{B(f,\gamma+\delta/4)\cup\F^{\delta/4}}\,.
\]
Indeed, since the points in $\F^{\delta/4}$ form a $\delta/4$ cover of $\F$,
all the points in $(B(f,\gamma)\cup\F)^{\delta/2}$ are at distance less than
$\delta/4$ of one and only one point of $\F^{\delta/4}$ (the unicity comes
from the fact that they are $\delta/2$ separated).
We can thus establish an injection from points in
$(B(f,\gamma)\cup\F)^{\delta/2}$ to corresponding points in $\F^{\delta/4}$
and the image of this injection is included in $B(f,\gamma+\delta/4)$ since
the image points are within distance $\delta/4$ of points in $B(f,\gamma)$.

Now we obtain
\[
N((B(f',\delta+\delta/4)\cup\F)^{\delta/2}, \delta/4)
\leq N(B(f',3\delta/2)\cup\F^{\delta/4}, \delta/8)\,.
\]
We thus have
\begin{eqnarray*}
N(B(f,\delta)\cup \F,\delta/2)
&\leq& N(B(f',\delta+\delta/4)\cup\F, \delta/2)\\
&\leq& N((B(f',\delta+\delta/4)\cup\F)^{\delta/2}, \delta/4)\\
&\leq& N(B(f',3\delta/2)\cup\F^{\delta/4}, \delta/8)\,.
\end{eqnarray*}
This gives
\begin{eqnarray*}
\sup_{f\in \F} N(B(f,\delta)\cap\F,\delta/2)&\leq& \sup_{f\in
\F^{\delta/4}} N(B(f,3\delta/2)\cap\F^{\delta/4},\delta/8)\\
&=& N(\F^{\delta/4},3\delta/2,\delta/8)\,.
\end{eqnarray*}
Hence
\[
H(\F,\delta,\delta/2)\leq H(\F^{\delta/4},3\delta/2,\delta/8)\,.
\]
By the same argument as in previous Lemma we obtain
\[
\frac{\delta}{8}H^{1/2}(\F^{\delta/4},3\delta/2,\delta/8) \leq
K\omega(\F^{\delta/4},3\delta/2)\,.
\]
\qed
\end{proof}

\subsection{Applications}
\begin{example}
If for all $\eps>0$,
\[
N(\F,\eps) \leq \eps^{-V}\,,
\]
then for all $\eps>0$,
\[
H(\conv(\F), \eps) \leq \eps^{-2V/(2+V)}\log^{2V/(2+V)}\eps^{-1}\,.
\]
\end{example}
\begin{proof}
Recall from Example \ref{ex:e1} that
\[
\omega(\conv(\F),\delta) \leq K\delta^{2/(2+V)}\log^{V/(2+V)} \delta^{-1}\,.
\]
Now, using Lemma \ref{le:lif} we get 
\begin{eqnarray*}
H^{1/2}(\conv(\F),2^{-k}) &\leq& K\sum_{i=0}^k 2^{i}
2^{2(1-i)/(2+V)}(i-1)^{V/(2+V)}\\ &=& K\sum_{i=0}^k
(2^{V/(2+V)})^i(i-1)^{V/(2+V)}\,.
\end{eqnarray*}
We check that in the above sum, the $i$-th term is always larger than
twice the $i-1$-th term (for $i\geq 2$) so that we can upper bound the
sum by the last term,
\[
H^{1/2}(\F,2^{-k}) \leq K(2^{V/(2+V)})^k(k-1)^{V/(2+V)}\,,
\]
hence, using $\eps=2^{-k}$, we get the result.\qed
\end{proof}
Note that the result we obtain contains an extra logarithmic factor compared
to the optimal bound \cite{vdVW96,mendelson01}.

\begin{example}
If for all $\eps>0$,
\[
H(\F,\eps) \leq \eps^{-V}\,,
\]
then for all $\eps>0$, for $0<V<2$,
\[
H(\conv(\F), \eps) \leq \eps^{-2}\log^{1-V/2}\eps^{-1}\,,
\]
for $V=2$,
\[
H(\conv(\F), \eps) \leq \eps^{-2}\log^2\eps^{-1}\,,
\]
and for $V>2$,
\[
H(\conv(\F), \eps) \leq \eps^{-V}\,.
\]
\end{example}
\begin{proof}
The proof is similar to the previous one.
\qed
\end{proof}

In this example, all the bounds are known to be sharp \cite{carl01,gao01}.


\begin{thebibliography}{10}
\bibitem{bp90}
K.~Ball and A.~Pajor.
\newblock The entropy of convex bodies with ``few'' extreme points.
\newblock {\em London MAth. Soc. Lectrure Note Ser.} 158, pages 25--32, 1990.

\bibitem{barboumen02}
P.~Bartlett, O.~Bousquet and S.~Mendelson.
\newblock Localized Rademacher Complexity.
\newblock {\em Preprint}, 2002.

\bibitem{boulugmas02}
S.~Boucheron, G.~Lugosi and P.~Massart.
\newblock Concentration inequalities using the entropy method.
\newblock {\em Preprint}, 2002.

\bibitem{bousquet02}
O.~Bousquet.
\newblock A Bennett concentration inequality and its application to
empirical processes.
\newblock {\em Comptes Rendus de l'Acad\'emie des Sciences}, 2002.

\bibitem{carl97}
B.~Carl.
\newblock Metric entropy of convex hulls in Hilbert spaces.
\newblock {\em Bulletin of the London Mathematical Society}, 29, pages
452--458, 1997.

\bibitem{carl01}
B.~Carl, I.~Kyrezi and A.~Pajor.
\newblock Metric entropy of convex hulls in Banach spaces.
\newblock {\em Journal of the London Mathematical Society}, 2001.

\bibitem{cre01}
J.~Creutzig and I.~Steinwart.
\newblock Metric entropy of convex hulls in type $p$ spaces -- the critical
case.
\newblock 2001.

\bibitem{dudley87}
R.~Dudley.
\newblock Universal Donsker classes and metric entropy.
\newblock {\em Annals of Probability}, 15, pages 1306--1326, 1987.

\bibitem{dudley00}
R.~Dudley.
\newblock Uniform central limit theorems.
\newblock Cambridge University Press, 2000.

\bibitem{dud}
R.~Dudley.
\newblock {\em Private communication}, 2001.

\bibitem{gao01}
F.~Gao.
\newblock Metric entropy of convex hulls.
\newblock {\em Israel Journal of Mathematics}, 123, pages 359--364, 2001.

\bibitem{kolpan00}
V.~I.~Koltchinskii and D.~Panchenko.
\newblock Rademacher processes and bounding the risk of function
learning.
\newblock In {\em High Dimensional Probability II}, Eds. E.Gine, D.Mason
and J.Wellner, pp. 443 - 459, 2000.

\bibitem{lt91}
M.~Ledoux and M.~Talagrand
\newblock Probability in Banach spaces.
\newblock Springer-Verlag, 1991.

\bibitem{lilinde01}
W.~Li and W.~Linde.
\newblock Metric entropy of convex hulls in Hilbert spaces.
\newblock {\em Preprint}, 2001.

\bibitem{lif}
M.~Lifshits.
\newblock Gaussian random functions.
\newblock Kluwer, 1995.

\bibitem{massart00}
P.~Massart.
\newblock Some applications of concentration inequalities to statistics.
\newblock {\em Annales de la Facult\'e des Sciences de Toulouse},
IX:245-303, 2000.

\bibitem{mendelson01}
S.~Mendelson.
\newblock On the size of convex hulls of small sets.
\newblock {\em Preprint}, 2001.

\bibitem{mendelson02}
S.~Mendelson.
\newblock Improving the sample complexity using global data.
\newblock {\em Preprint}, 2001.

\bibitem{vdVW96}
A.~van der Vaart and J.~Wellner.
\newblock Weak convergence and empirical processes with applications to
statistics.
\newblock John Wiley \& Sons, New York, 1996.
\end{thebibliography}
\end{document}